\def\UseSection{
      \numberwithin{equation}{section}
	\theoremstyle{plain}
      \newtheorem{theorem}    {Theorem}[section]
      \DefineTheorems 
}

\def\DefineTheorems{
	
	\newtheorem{lemma}      [theorem] {Lemma}
	
	\newtheorem{prop}       [theorem] {Proposition}
	
	\newtheorem{cor}        [theorem] {Corollary}

	\theoremstyle{definition}
	\newtheorem{defn}       [theorem] {Definition}

	\theoremstyle{definition}

}

\newcommand{\bt}   {\begin{theorem}}
\newcommand{\et}   {\end  {theorem}}
\newcommand{\bl}   {\begin{lemma}}
\newcommand{\el}   {\end  {lemma}}
\newcommand{\bp}   {\begin{prop}}
\newcommand{\ep}   {\end  {prop}}
\newcommand{\bc}   {\begin{cor}}
\newcommand{\ec}   {\end  {cor}}
\newcommand{\bd}   {\begin{defn}}
\newcommand{\ed}   {\end  {defn}}

\newcommand{\ba}   {\begin{array}}
\newcommand{\ea}   {\end  {array}}
\newcommand{\be}   {\begin{enumerate}}
\newcommand{\ee}   {\end  {enumerate}}
\newcommand{\bi}   {\begin{itemize}}
\newcommand{\ei}   {\end  {itemize}}

\def\eq#1\en{\begin{equation}#1\end{equation}}  
\def\eqsplit#1\ensplit{
	\begin{equation}\begin{split}#1\end{split}\end{equation}
	}
\def\eqalign#1\enalign{
	\begin{align}#1\end{align}
	}
\def\eqmul#1\enmul{
	\begin{multline}#1\end{multline}
	}
\newcommand{\eqarrstar} {\begin{eqnarray*}} 
\newcommand{\enarrstar} {\end{eqnarray*}} 
\newcommand{\eqarray}   {\begin{eqnarray}} 
\newcommand{\enarray}   {\end{eqnarray}} 
 
\makeatletter
\newcommand{\labelcounter}[2]{{%
	\stepcounter{#1}
	\protected@write\@auxout{}%
	{\string\newlabel{#2}{{\csname the#1\endcsname}{\thepage}}}%
	{\ref{#2}}
	}}
\makeatother
\newcommand{\Cbold} {{\mathbb C}}

\newcommand{\Nbold} {{\mathbb N}}

\newcommand{\Zbold} {{\mathbb Z}}

\newcommand{\spose}[1] {{\hbox to 0pt{#1\hss}} }
\newcommand{\ltapprox} {\mathrel{\spose{\lower 3pt\hbox{$\mathchar"218$}}
\raise 2.0pt\hbox{$\mathchar"13C$}}}
\newcommand{\gtapprox} {\mathrel{\spose{\lower 3pt\hbox{$\mathchar"218$}}
\raise 2.0pt\hbox{$\mathchar"13E$}}}

\documentclass[12pt]{amsart}
\usepackage{enumitem}
\usepackage{amsmath,amssymb,amsthm}
\usepackage[dvips]{graphicx}
\usepackage{eepic}
\newtheorem{THM}{Theorem}[section]
\newtheorem{COR}[THM]{Corollary}

\newtheorem{MOD}[THM]{Model}

\newtheorem{LEM}[THM]{Lemma}
\newtheorem{PRP}[THM]{Proposition}
\newtheorem{DEF}[THM]{Definition}

\UseSection  
\newcommand{\hlf}{\frac{1}{2}}
\newcommand{\ra}{\rightarrow}
\newcommand{\lra}{\leftrightarrow}

\renewcommand{\to}      {\rightarrow}

\newcommand{\del}{\partial}

\newcounter{countC}  
\newcounter{countR}  
\newcommand{\re}{\mathbb{R}}
\newcommand{\Z}{\Zbold}
\newcommand{\N}{\Nbold}
\newcommand{\C}{\Cbold}

\newcommand{\mc}[1]{\mathcal{#1}}
\newcommand{\mP}{\mathbb{P}}

\newcommand{\mE}{\mathbb{E}}

\newcommand{\mG}{\mc{G}}

\newcommand{\UD}{\updownarrow}
\newcommand{\LR}{\leftrightarrow}
\newcommand{\WE}{\LR}
\newcommand{\NS}{\UD}

\newcommand{\smallE}{\scriptstyle \rightarrow}
\newcommand{\smallW}{\scriptstyle \leftarrow}

\newcommand{\ssmallW}{\scriptscriptstyle \leftarrow}
\newcommand{\smallN}{\scriptstyle \uparrow}
\newcommand{\smallS}{\scriptstyle \downarrow}

\newcommand{\smallnw}{\scriptscriptstyle \nwarrow}
\newcommand{\smallsw}{\scriptscriptstyle \swarrow}

\newcommand{\NE}{\begin{picture}(,)
\put(2,-5){$\rightarrow$}
\put(0,.5){$\uparrow$}
\end{picture}\hspace{.5cm}
}
\newcommand{\smallNE}{\begin{picture}(,)
\put(1.5,-3){$\smallE$}
\put(0,.5){$\smallN$}
\end{picture}\hspace{.35cm}
}

\newcommand{\smallSE}{\begin{picture}(,)
\put(1.4,3.2){$\smallE$}
\put(,-.5){$\smallS$}
\end{picture}\hspace{.35cm}
}
\newcommand{\SW}{\begin{picture}(,)
\put(0,4.8){$\leftarrow$}
\put(8.2,-0.5){$\downarrow$}
\end{picture}\hspace{.5cm}
}
\newcommand{\smallSW}{\begin{picture}(,)
\put(0.5,3.2){$\smallW$}
\put(6.3,-0.5){$\smallS$}
\end{picture}\hspace{.35cm}
}
\newcommand{\ssmallSW}{\hspace{-.4cm}\begin{picture}(,)
\put(0.5,3.2){$\smallW$}
\put(6.3,-0.5){$\smallS$}
\end{picture}\hspace{.35cm}
}

\newcommand{\smallNW}{\begin{picture}(,)
\put(0.5,-3){$\smallW$}
\put(6.3,.5){$\smallN$}
\end{picture}\hspace{.35cm}
}

\newcommand{\SWE}{\begin{picture}(,)
\put(0,5){$\leftarrow$}
\put(5,5){$\rightarrow$}
\put(5.5,-0.5){$\downarrow$}
\end{picture}\hspace{.6cm}
}

\newcommand{\smallOTSP}{\begin{picture}(,)
\put(6.2,0){$\smallN$}
\put(.3,-3.3){$\smallW$}
\put(1.7,0){$\smallnw$}
\end{picture}\hspace{.4cm}
}
\newcommand{\OTSP}{\begin{picture}(,)
\put(0.5,-5){$\leftarrow$}
\put(0.5,.2){$\nwarrow$}
\put(8.7,0.5){$\uparrow$}
\end{picture}\hspace{.5cm}
}

\newcommand{\smallFSOSP}{\begin{picture}(,)
\put(5.5,2.7){$\smallN$}
\put(5.5,-2.3){$\smallS$}
\put(.9,3.1){$\smallnw$}
\put(.9,-2){$\smallsw$}
\put(0.5,.5){$\ssmallW$}
\end{picture}\hspace{.35cm}
}

\newcommand{\blank}[1]{}

\newcommand{\barc}{\bar{\mathcal{C}}}

\newcommand{\Qed}{\qed \bigskip}

\newcommand{\nra}{\nrightarrow}
\usepackage[usenames]{color}

\begin{document}

\title  {
      Forward clusters for degenerate random environments}

\author[Holmes]{Mark Holmes}
\address{Department of Statistics, University of Auckland}
\email{mholmes@stat.auckland.ac.nz}
\author[Salisbury]{Thomas S. Salisbury} 
\address{Department of Mathematics and Statistics, York University}
\email{salt@yorku.ca}


\begin{abstract}
We consider connectivity properties and asymptotic slopes for certain random directed graphs on $\Z^2$ in which  the set of points $\mc{C}_o$ that the origin connects to is always infinite.  We obtain conditions under which the complement of $\mc{C}_o$ has no infinite connected component. 
Applying these results to one of the most interesting such models leads to an improved lower bound for the critical occupation probability for oriented site percolation on the triangular lattice in 2 dimensions. 
\end{abstract}

\maketitle

\section{Introduction and Main Results}
\label{sec:intro}
The main objects of study in this paper are the {\em 2-dimensional orthant model} (one of the most interesting examples within a class of models called {\em degenerate random environments}), and its dual model, a version of {\em oriented site percolation}.  Part of the motivation for studying degenerate random environments is an interest in the behaviour of random walks in random environments that are non-elliptic.  
Indeed, many of the results of this paper and of \cite{HS_DRE1} have immediate implications for the behaviour (in particular, directional transience) of random walks in certain non-elliptic environments (see~e.g.~\cite{HS_RWDRE}).  

For fixed $d\ge 2$, let $\mc{E}=\{\pm e_i: i=1,\dots,d\}$ be the set of unit vectors in $\Z^d$, and let $\mc{P}$
denote the power set of $\mc{E}$.
Let $\mu$ be a probability measure on $\mc{P}$.  A {\em degenerate random environment} (DRE) is a random directed graph, 
i.e.~an element $\mc{G}=\{\mc{G}_x\}_{x\in \Z^d}$ of $\mc{P}^{\Z^d}$.  We equip $\mc{P}^{\Z^d}$ with the product $\sigma$-algebra and the product measure $\mP=\mu^{\otimes \Z^d}$, so that $\{\mc{G}_x\}_{x\in \Z^d}$ are i.i.d.~under $\mP$.  We denote the expectation of a random variable $Z$ with respect to $\mP$ by $\mE[Z]$. 

We say that the DRE is {\em $2$-valued} when $\mu$ charges exactly two points, i.e.~there exist distinct $E_1,E_2\in\mc{P}$ and $p \in (0,1)$ such that $\mu(\{E_1\})=p$ and $\mu(\{E_2\})=1-p$.  As in the percolation setting, there is a natural coupling of graphs for all values of $p$ as follows.  Let $\{U_x\}_{x \in \Z^d}$ be i.i.d.~standard uniform random variables under $\mP$.  Setting
\begin{equation}
\label{eq:coupling}
\mc{G}_x=\begin{cases}
E_1, &\text{ if }U_x<p,\\
E_2, &\text{otherwise,}
\end{cases}
\end{equation}
it is easy to see that for any fixed $p$, $\mc{G}$ has the correct law, and that the set of $E_1$ sites is increasing in $p$, $\mP$-almost surely.

For the most part, in this paper, we will consider 2-dimensional models.  Interesting examples of $2$-dimensional, $2$-valued models include the following:

\begin{MOD} $(\NE\SW)$:
\label{exa:NE_SW}
Let $E_1=\{\uparrow,\rightarrow\}$ and $E_2=\{\downarrow,\leftarrow\}$ (and set $\mu(\{E_1\})=p$, $\mu(\{E_2\})=1-p$). 
\end{MOD}
We call the generalization to $d$ dimensions the {\it orthant model} (so this is the {\it 2-$d$ orthant model}).
\begin{MOD} $(\SWE\uparrow)$
\label{exa:SWE_N}
Let $E_1=\{\leftarrow,\downarrow,\rightarrow\}$ and $E_2=\{\uparrow\}$.
\end{MOD}
\begin{MOD} $(\NS \WE)$:
\label{exa:NS_WE}
Let $E_1=\{\uparrow,\downarrow\}$ and $E_2=\{\rightarrow,\leftarrow\}$.
\end{MOD}
We will focus on Model \ref{exa:NE_SW} in this paper, but we believe extensions of our results to Model \ref{exa:SWE_N} are possible. A central limit theorem for random walk in the degenerate random environment of Model \ref{exa:NS_WE} is proved in \cite{BD11}. See \cite{HS_DRE1} for more on these and other examples of $2$-dimensional $2$-valued models, and their properties.  For any set $A$, let $|A|$ denote its cardinality.

\begin{DEF}
\label{def:connections}
Given an environment $\mc{G}$:
\begin{itemize}
\item We say that $x$ is {\em connected} to $y$, and write $x\ra y$ if: there exists an $n\ge 0$ and a sequence $x=x_0,x_1,\dots, x_n=y$ such that $x_{i+1}-x_{i}\in \mc{G}_{x_i}$ for $i=0,\dots,n-1$. We say that $x$ and $y$ are {\em mutually connected}, or that they {\em communicate}, and write $x\lra y$ if $x\ra y$ and $y \ra x$.  
\item  Define $\mc{C}_x=\{y\in \Z^d:x \ra y\}$ (the {\em forward} cluster), $\mc{B}_y=\{x\in\Z^d:x\ra y\}$ (the {\em backward} cluster), and $\mc{M}_x=\{y \in \Z^d:x\lra y\}=\mc{B}_x\cap\mc{C}_x$ (the {\em bi-connected} cluster). \newline
Set
$\theta_+=\mP(|\mc{C}_o|=\infty)$, $\theta_-=\mP(|\mc{B}_o|=\infty)$, and $\theta=\mP(|\mc{M}_o|=\infty)$.
\item A nearest neighbour path in $\Z^d$ is {\em open in $\mc{G}$} if that path consists of directed edges in $\mc{G}$.
\end{itemize}
\end{DEF}
It is interesting to consider percolation-type questions where $\theta_+<1$, and indeed in the true percolation settings where $\mu(\{\emptyset\})=1-\mu(\{E\})$ for some configuration $E$, there is a simple relation between the directed percolation probabilities of the form $\theta_+=\mu(\{E\})\theta_-$ (see \cite[Lemma 2.5]{HS_DRE1}).  Our current interest lies in those cases where $\mu$ is such that 
\begin{equation}
\boxed{\theta_+=1}.
\end{equation}
This is precisely the condition which ensures that a random walk on this random graph visits infinitely many sites. It is equivalent to the condition that there exists a set of orthogonal directions $V\subset \mc{E}$ such that $\mu(\mc{G}_o\cap V\ne \varnothing)=1$, i.e.~that almost surely at every site, some element of $V$ occurs \cite[Lemma 2.2]{HS_RWDRE}.  Models \ref{exa:NE_SW}-\ref{exa:NS_WE} above satisfy this condition (e.g.~by taking $V=\{\uparrow,\leftarrow\}$).  Note that $\theta_-=1$ if and only if there exists $e \in \mc{E}$ such that $\mu(e \in \mc{G}_o)=1$, and that none of these three models satisfy this condition.  In fact $\theta_-\in (0,1)$ for each of these models when $p \in (0,1)$ \cite{HS_DRE1}.  

Model \ref{exa:NE_SW} exhibits a phase transition (in fact two phase transitions by symmetry) when the parameter $p$ changes \cite{HS_DRE1}.  While $\theta_+=1$ and $\theta_->0$ for all $p$, the geometry of an infinite $\mc{B}_x$ changes from having a non-trivial boundary to being all of $\Z^2$ as $p$ decreases from $1$ to $\hlf$.  Moreover the critical point $p_c$ at which this transition takes place is also the critical point $p_c^{\smallOTSP}$ for oriented site percolation on the triangular lattice (OTSP).  Information about the geometry of an infinite $\mc{M}_x$ is then inferred, based on the geometry of $\mc{B}_x$, and crude estimates of $p_c$ are given.  Similar results hold for Model  \ref{exa:SWE_N}, which is dual to a partially oriented site percolation model on the triangular lattice.  
However, for Model \ref{exa:NS_WE}, 
if $|\mc{B}_x|=\infty$ then $\mc{B}_x=\Z^2$ almost surely, regardless of $p$ \cite{HS_DRE1}.

Simulations indicate that $\mc{C}_o$ and infinite $\mc{B}_o$ clusters have similar geometry, except that $\mc{C}_o$ typically has ``holes" whereas $\mc{B}_o$ does not.  In order to give a clearer description of this weak kind of duality, we study the geometry of
$\barc_x\supset \mc{C}_x$, defined by
\begin{multline}
\label{eq:Cbar}
\barc_x=\{z \in \Z^d: \text{ every infinite nearest-neighbour self-avoiding} \\  \text{path starting at $z$ passes through $\mc{C}_x$}\}.
\end{multline}

\begin{figure}
\begin{center}
\includegraphics[scale=.45]{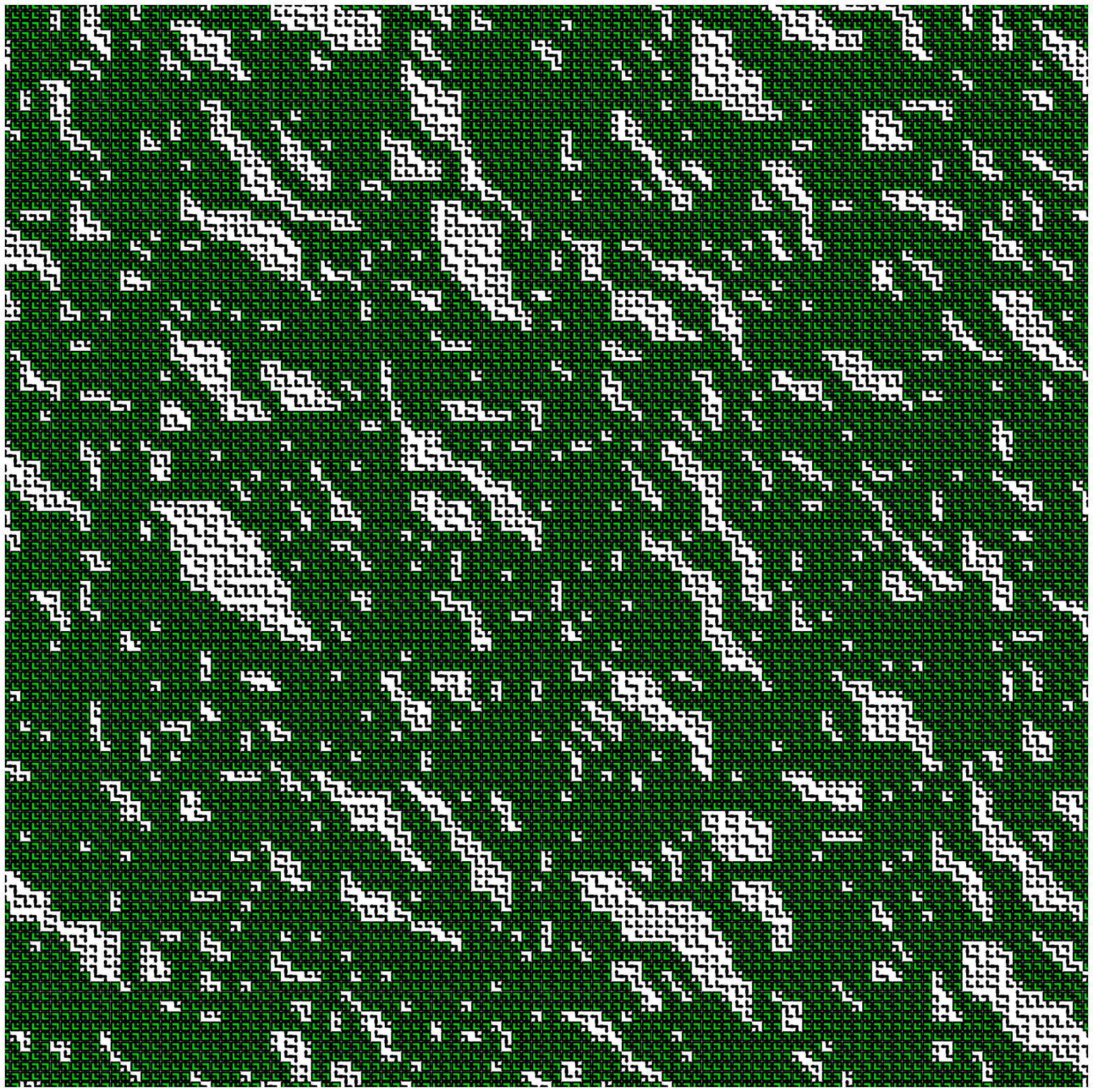} 
\hspace{.2cm}
\includegraphics[scale=.45]{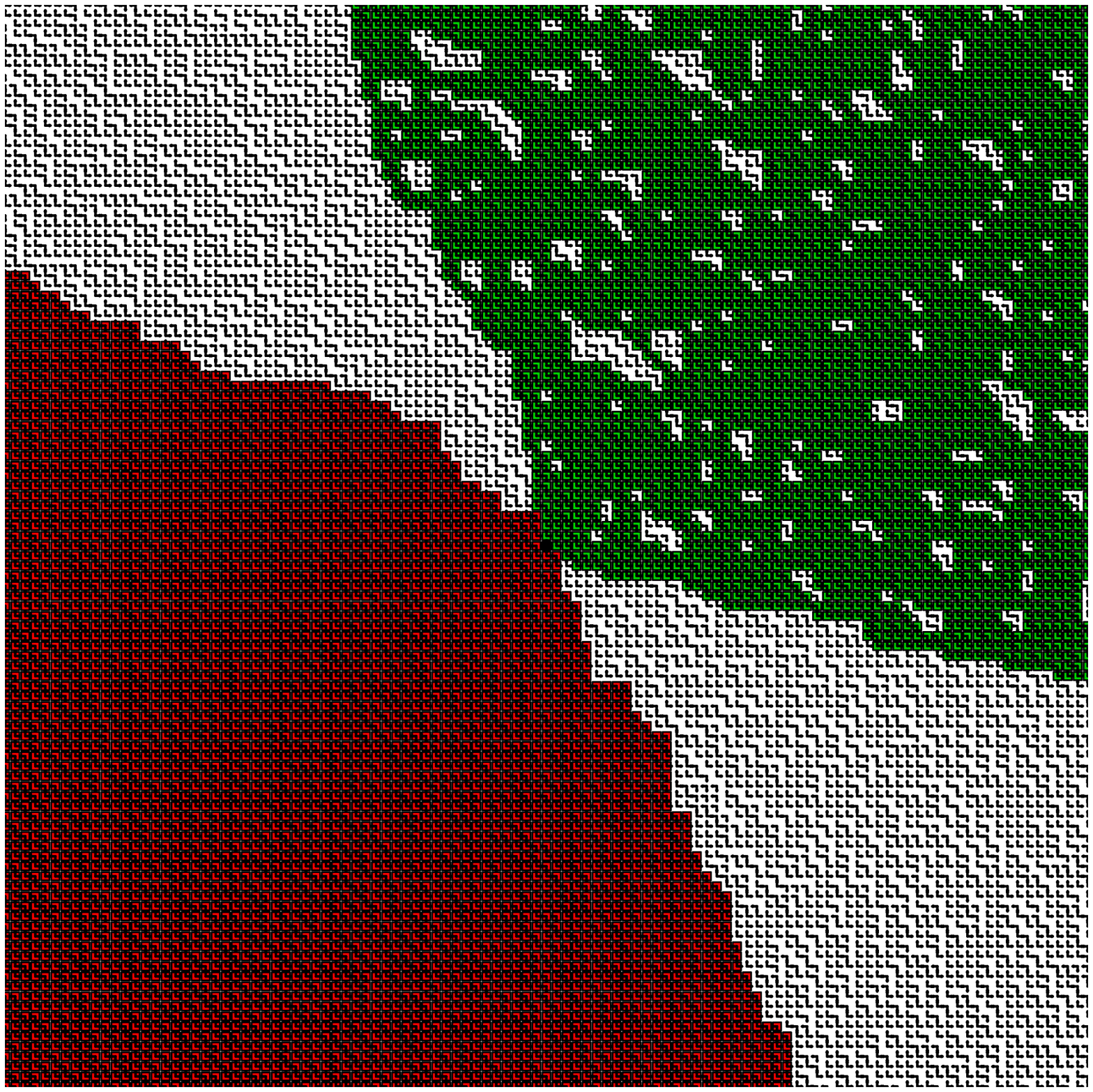} 
\end{center}
\caption{Simulations of [$\mc{C}_o$ when $p=.5$] and of [$\mc{C}_o$ and $\mc{B}_o$ when $p=.7$] for Model \ref{exa:NE_SW}.}
\label{fig:orthant_C_o_B_o}
\end{figure}

The results of this paper can be separated into two groups.  The first group concerns the possible geometries of the sets $\barc_x$ for a class of 2-dimensional models.  These appear in Section \ref{sec:Cgeneral}, and include adaptations of the results about the geometries of $\mc{B}_x$ clusters in \cite{HS_DRE1}.  

The second group of results concerns Model \ref{exa:NE_SW}.    
Some of these results, included in Section \ref{sec:orthant_main}, are again adaptations of results in \cite{HS_DRE1}, and are proved using similar arguments.  
Others exploit much more deeply the duality between $\mc{B}_o$, $\barc_o$ and OTSP.  An analysis of OTSP (see Section \ref{sec:OP} and also e.g.~\cite{Dur84}) is used to prove the following result (see also Figure \ref{fig:orthant_C_o_B_o}) in Sections \ref{sec:orthant_main} and \ref{sec:slopes}.
\begin{THM}
\label{thm:orthant_main}
For Model \ref{exa:NE_SW} $(\NE\SW)$, the following hold:
\begin{enumerate}[label={\emph{(\Roman*)}},ref={\theTHM.(\Roman*)}]
\item $1-p_c\le p\le p_c$ if and only if $\mc{M}_o$ is infinite with positive probability. \label{thm:Mphasetransition}
\item $\barc_o(p_1)\supsetneq \barc_o(p_2)$ $\mP$-a.s. whenever
$p_c\le p_1<p_2\le 1$, under the coupling (\ref{eq:coupling}). \label{thm:orthant_monotone}
\item when $p>p_c$, the northwest-pointing boundary of $\barc_o$ has an asymptotic slope $\rho_p< -1$, with $\rho_p$ strictly increasing to $-1$ as $p\downarrow p_c$. The southeast-pointing boundary of $\barc_o$ has asymptotic slope $1/\rho_p$. 
\label{thm:orthant_slope}
\end{enumerate}
\end{THM}
Note that results of \cite{HS_DRE1} already show $\Rightarrow$ for (I).
But $\Leftarrow$ is only verified there assuming a stronger condition on $\mc{M}_o$ (that $\mc{M}_o$ be ``gigantic"). 
In contrast to (II), the cluster $\mc{C}_o(p)$ is {\em not} monotone in $p$ under the natural coupling (\ref{eq:coupling}).  For example, letting  $\phi_x(p)=\mP_p\big(x\in \mc{C}_o\big)$, we see that $\phi_{(-1,1)}(p)=0$ if $p=0$ or $p=1$, and $\phi_{(-1,1)}(p)>0$ otherwise (in fact $\phi_{(-1,1)}(1/2)>1/2$).  Even when restricting attention to $p\ge 1/2$, holes may open or close in $\mc{C}_o$ as $p$ increases.  Finally note that a corollary of (III) is that random walks in random environments whose supports are $\NE$ with probability $p$ and $\SW$ with probability $1-p$ are transient in direction $(1,1)$ when $p>p_c$ (see \cite{HS_RWDRE}).  
On the other hand, by finding arbitrarily large connected circuits in Model \ref{exa:NE_SW} when $0.427\le p\le 0.573$ we obtain the following result (proved in Section \ref{sec:pcbounds}) about OTSP.
\begin{THM}
\label{thm:OTSP_pc}
The critical occupation probability for oriented site percolation on the triangular lattice ($d=2$) is at least $0.5730$.
\end{THM}
This improves on the best rigorous bounds that we have found in the literature, namely:  $0.5466\le p_c^{\smallOTSP}\le 0.7491$ (see \cite{HS_DRE1,BBS}).  Note that the estimated value is $p_c\approx 0.5956$ \cite{DeBE,JG}.   We believe that an adaptation of these arguments to Model \ref{exa:SWE_N} yields a bound on the critical occupation probability $p_c^{\smallFSOSP}$ for a partially oriented site percolation model, and it may be that one can similarly obtain useful bounds on the critical probability $p_c^{\smallNW}$ for oriented site percolation on the square lattice via related degenerate random environments on the triangular lattice.

\section{General models: The forward cluster $\mc{C}_o$}
\label{sec:Cgeneral}
In this section we investigate properties of the random sets $\mc{C}_x$.  Typically these clusters are rather different from connected clusters in percolation models (where $E_2=\varnothing$).  For example, in the $2$-valued setting the sets  $\mc{C}_y$ are not increasing under the coupling (\ref{eq:coupling}) unless $E_1\supset E_2$.  In particular for Models \ref{exa:NE_SW}-\ref{exa:NS_WE} the cluster $\mc{C}_o$ is not monotone in $p$.  
A natural question is to ask whether or not the connection events $\{x\in \mc{C}_o\}_{x\in \Z^d}$ are positively correlated, i.e.~whether $\mP(y \in \mc{C}_o,x\in \mc{C}_o)\ge \mP(y \in \mc{C}_o)\mP(x \in \mc{C}_o)$.  While such a property is true in the percolation setting, this fails in general for degenerate random environments.  The easiest example is the case $E_1=\{e_1\}$ and $E_2=\{e_2\}$, where $\mP(e_1 \in \mc{C}_o,e_2\in \mc{C}_o)=0\ne p(1-p)=\mP(e_1 \in \mc{C}_o)\mP(e_2\in \mc{C}_o)$.  We believe that it fails for Models \ref{exa:NE_SW}-\ref{exa:NS_WE} (and many others) as well.

On the other hand, for a fixed $y$, the events $\{y\in \mc{C}_x\}_{x \in \Z^d}$ are positively correlated. 
\begin{LEM}
\label{lem:jointC}
For $x,y\in \Z^d$, 
\eqalign
\mP(y \in \mc{C}_o,y \in \mc{C}_x)\ge &\mP(y \in \mc{C}_o)\mP(y \in \mc{C}_x).
\enalign
\end{LEM}
\proof  
It is sufficient to prove that
$P(y \in \mc{C}_o,y \notin \mc{C}_x)\le \mP(y \in \mc{C}_o)\mP(y \notin \mc{C}_x)$.
Let $B(n)$ be the set of lattice sites in a ball of radius $n$, centred at $o$.  Let $E_n=\{o \ra y \text{ in }B(n)\}$ be the event that there is a path from $o$ to $y$ lying entirely in $B(n)$.  Let $F_n=\{x\rightarrow y \text{ in }B(n)\}^c$ be the event that there is no path connecting $x$ to $y$ that lies entirely in $B(n)$.  Then $E_n$ and $F_n$ are increasing and decreasing events respectively, with $E\equiv \{o\ra y\}=\bigcup_{n=1}^{\infty}E_n$ and $F\equiv \{x\nra y\}=\bigcap_{n=1}^{\infty}F_n$.  Therefore,
\begin{align*}
\mP(E\cap F)=\lim_{M\ra \infty}\mP(E_M\cap F)\le \liminf_{M\ra \infty}\mP(E_M\cap F_M).
\end{align*}
Let $\mc{C}_{x,M}$ be the set of sites that can be reached from $x$ using only sites in $B(M)$.  Observe that if  $z\in \mc{C}_{x,M}$ and $z\ra y$ in $B(M)$ then $x\ra y$ in $B(M)$.  Thus for any $C\subset B(M)$ with $y\notin C$, on the event $\{\mc{C}_{x,M}=C\}$, we have that $\{o \ra y \text{ in }B(M)\}$ occurs if and only if $\{o \ra y \text{ in }B(M)\setminus C\}$ occurs.  This latter event depends only on the random variables $\{\mc{G}_z:z\in B(M)\setminus C\}$, while $\{\mc{C}_{x,M}=C\}$ depends only on the random variables $\{\mc{G}_z:z\in C\}$.  Thus we have
\begin{align*}
\mP(E_M\cap F_M)=&\sum_{C\subset B(M)\text{ s.t. } y\notin C}\mP\big(0 \ra y \text{ in }B(M),\mc{C}_{x,M}=C\big)\\
=&\sum_{C\subset B(M)\text{ s.t. } y\notin C}\mP\big(0 \ra y \text{ in }B(M)\setminus C,\mc{C}_{x,M}=C\big)\\
=&\sum_{C\subset B(M)\text{ s.t. } y\notin C}\mP\big(0 \ra y \text{ in } B(M)\setminus C\big)\mP\big(\mc{C}_{x,M}=C\big)\\
\le &\sum_{C\subset B(M)\text{ s.t. } y\notin C}\mP\big(0 \ra y\big)\mP\big(\mc{C}_{x,M}=C\big)=\mP\big(0 \ra y\big)\mP\big(x\nra y \text{ in }B(M)\big).
\end{align*}
In other words, $\mP(E_M\cap F_M)\le \mP\big(0 \ra y\big)\mP(F_M)$, and taking the limit as $M\ra \infty$ establishes the result.
 \Qed

Note that by translation invariance and relabelling of vertices, Lemma \ref{lem:jointC} is equivalent to saying that $\mP(x \in \mc{B}_o,y\in \mc{B}_o)\ge \mP(x \in \mc{B}_o)\mP(y \in \mc{B}_o)$
so roughly speaking, knowing that something connects to $o$ makes it more likely that other things connect to $o$.    

Let $C\subset\Z^2$. We say that $C$ has a {\em finite hole} $G$ if $G\subset \Z^2\setminus C$, $G$ is finite, $G$ is connected in $\Z^2$, and every $z\in\Z^2\setminus G$ that is a neighbour of $G$ must belong to $C$. 
The following elementary lemma implies that $\barc_o$ is obtained from $\mc{C}_o$ by filling in all finite holes, and that the backward cluster of a finite hole $G$ is simply $G$.
\begin{LEM}
\label{lem:Choles}
Suppose that $x\in \barc_o\setminus \mc{C}_o$.  Then $x$ belongs to a finite hole $G$ of $\mc{C}_o$ and $\mc{B}_x\subset G$.
\end{LEM}
\proof Let $x\in \barc_o\setminus \mc{C}_o$, and let $G\subset\Z^2\setminus \mc{C}_o$ be the $\Z^2$-connected cluster of $x$ in $\Z^2\setminus \mc{C}_o$.  Clearly every neighbour $z\in\Z^2\setminus G$ of $G$ is in $\mc{C}_0$.   Let $y \in G$.  Then there is a finite self-avoiding path in $G$ from $x\in \barc_o$ to $y$, so $y \in \barc_o$.  Thus $G\subset \barc_o\setminus \mc{C}_o$.  Every infinite connected subset of $\Z^2$ contains an infinite nearest-neighbour self-avoiding path (it is easy to construct this iteratively - if $x_0$ is connected to infinity then there exists some neighbour $x_1$ of $x_0$ that connects to infinity off $\{x_0\}$ etc). Therefore $G$ is finite.  Thus $G$ is a finite hole containing $x$.
Finally, since $x\in\mc{B}_x$ and $\mc{B}_x\cap \mc{C}_o=\emptyset$, it follows that $\mc{B}_x\subset G$. 
\Qed

The above relationships between $\mc{B}$ and $\mc{C}$ clusters are rather weak.  We can however prove results about the $\mc{C}$ clusters which are dual to results about $\mc{B}$ clusters in \cite{HS_DRE1}, using modifications of the arguments from \cite{HS_DRE1}.  To state these results, we need the notion of blocking functions.

\begin{DEF}
\label{def:blocked}
Given $w:\Z\ra \Z$, define $w_{\le}\subset \Z^2$ and $w_{>}\subset \Z^2$ by
$$
w_\le = \{y\in \Z^2: y^{[2]}\le w(y^{[1]})\} \quad \text{ and }\quad w_> = \{y\in \Z^2: y^{[2]}> w(y^{[1]})\}.
$$
We say that $y$ is  \emph{below $w$} if $y^{[2]}\le w(y^{[1]})$, and \emph{strictly below $w$} if $y^{[2]}< w(y^{[1]})$. 
\newline We define $w_{\ge}$ and $w_<$ similarly, and speak likewise of $y$ being \emph{above $w$} or \emph{strictly above $w$}. 
\end{DEF}
We say that $w:\Z\ra \Z$ is a {\em forward lower blocking function} (flbf) for $\mc{G}$ if there is no open path in $\mc{G}$ from $w_{\ge}$ to $w_<$, i.e.~if $w_{\ge}\nrightarrow w_{<}$.  Similarly, $w$ is a {\em forward upper blocking function} (fubf) for $\mc{G}$ if $w_{\le}\nrightarrow w_{>}$.
Note that these notions are different from the (backward) lower blocking function $w_{<}\nrightarrow w_{\ge}$ and (backward) upper blocking function $w_{>}\nrightarrow w_{\le}$ in \cite{HS_DRE1}.  In particular, $w$ is a flbf if and only if $w-1$ is a bubf, and $w$ is a fubf if and only if $w+1$ is a blbf.

We write $\mc{A}_{E}=\{\mc{G}_o\cap E\ne \varnothing\}$ for $E\subset \mc{E}$ and use shorthand such as $\mc{A}_{\smallNW}=\mc{A}_{\{\smallW, \smallN\}}$.  The following Propositions are $\mc{C}$-dual versions of the $\mc{B}$ results \cite[Proposition 3.8 and Corollary 3.10]{HS_DRE1}.  
\begin{PRP}
\label{prp:Ctrichotomy}
Fix $d=2$.  Assume that $\mu(\mc{A}_{\smallNW})=1$,  $\mu(\mc{A}_{\smallSE})=1$, $\mu(\mc{A}_{\smallW})>0$, $\mu(\mc{A}_{\smallE})>0$, $\mu(\mc{A}_{\smallN})>0$, and $\mu(\mc{A}_{\smallS})>0$.  
\begin{enumerate}
\item 
The following $\mP$-a.s.~exhaust the possibilities for $\barc_x$:\smallskip
  \begin{enumerate}
  \item[(i)] $\barc_x=\Z^2$;
  \item[(ii)] There exists a decreasing flbf $W:\Z\to\Z$ such that $\barc_x=W_{\ge}$.
  \item[(iii)] There exists a decreasing fubf $W:\Z\to\Z$ such that $\barc_x=W_{\le}$;
  \end{enumerate}
\item Only one of (i), (ii), (iii) can have probability different from 0. 
\end{enumerate}
\end{PRP}
\begin{PRP}
\label{prp:Ctrichotomy2}
Fix $d=2$.  Assume that
$\mu(\mc{A}_{\smallNE})=1$,  $\mu(\mc{A}_{\smallNW})=1$, $\mu(\mc{A}_{\smallW})>0$, $\mu(\mc{A}_{\smallE})>0$, and $\mu(\mc{A}_{\smallN})>0$.  
\begin{enumerate}
\item 
The following $\mP$-a.s. exhaust the possibilities for $\barc_x$:\smallskip
  \begin{enumerate}
  \item[(i)] $\barc_x=\Z^2$;
  \item[(ii)] There exists a flbf $W:\Z\to\Z$ such that $\barc_x=W_{\ge}$;
  \end{enumerate}
\item Only one of (i) or (ii) can have probability different from 0. 
\end{enumerate}
\end{PRP}
Proposition \ref{prp:Ctrichotomy} applies to Model \ref{exa:NE_SW}, while Proposition \ref{prp:Ctrichotomy2} applies to Model \ref{exa:SWE_N}. 
When $\barc_x=W_{\ge}$ (so $W$ is a flbf) we say that $\barc_x$ is {\it blocked below}.  Similarly when $\barc_x=W_{\le}$ (so $W$ is a fubf) we say that $\barc_x$ is {\it blocked above}. \bigskip

An important notion that arises in the proofs of these results (and elsewhere throughout this paper) is the asymptotic slope of a path.  
\begin{DEF}
\label{def:slope}
A nearest-neighbour 
path $x_0,x_1,\dots$ with $x_i=(x_i^{[1]},x_i^{[2]})\in \Z^2$ is said to have \emph{asymptotic slope $\sigma$} if 
\[\lim_{n\ra \infty}\frac{x_n^{[2]}}{x_n^{[1]}}=\sigma.\]
\end{DEF}

\bigskip

\noindent {\em Proof of Proposition \ref{prp:Ctrichotomy}.}
  Take $x=o$. As in the proof of \cite[Proposition 3.8]{HS_DRE1}, we may construct NW or SE paths from any point. Suppose that $w\notin\mc{C}_o$, $\bar y\in\barc_o$, $\bar z\in\barc_o$,  $w^{[1]}=\bar y^{[1]}=\bar z^{[1]}$, but $\bar y^{[2]}<w^{[2]}<\bar z^{[2]}$. Because $\bar y$ is either in $\mc{C}_o$ or it is enclosed by $\mc{C}_o$, we can find $y\in\mc{C}_o$ such that $y^{[1]}=\bar y^{[1]}$ but $y^{[2]}\le \bar y^{[2]}$. Likewise there is a $z\in\mc{C}_o$ with $z^{[1]}=\bar z^{[1]}$ and $z^{[2]}\ge \bar z^{[2]}$. The SE paths from $y$ and $z$ intersect, by \cite[Lemma 2.3]{HS_DRE1}. So do the NW paths from $y$ and $z$. These four paths enclose $w$, so $w\in\barc_o$. Letting $L_i=\inf\{j:(i,j)\in\barc_o\}$ and $U_i=\sup\{j:(i,j)\in\barc_o\}$ it follows that $\barc_o=\{(i,j):L_i\le j\le U_i\}$. 

Case (i) corresponds to $L\equiv-\infty$ and $U\equiv \infty$. We can rule out the possibility that $L$ or $U$ jump from finite to infinite values, or vice versa, just as in \cite[Proposition 3.8]{HS_DRE1}. To see that if $L$ takes finite values, it must be decreasing, consider $\mc{G}_{(i,L_i)}$. By definition, $\downarrow\notin\mc{G}_{(i,L_i)}$. Since $\mu(\mc{A}_{\smallSE})=1$, it follows that $\rightarrow\in\mc{G}_{(i,L_i)}$. This implies that $L_{i+1}\le L_i$.  Similarly, $U_i$ is decreasing if finite. 

If $-\infty<L$ then $W=L$ is a flbf, and similarly if $U<\infty$ then $W=U$ is a fubf.  Thus it remains to prove that one of $L$ or $U$ must be infinite, and it suffices to do this for $L_0$ and $U_0$. To do this, we make use of a number of paths.  Define the $Nw$ path ($\subset \mc{C}_x$) from $x$ to be that path starting from $x$ obtained by following $\uparrow$ whenever possible, and otherwise following $\leftarrow$.  On a set of $\mP$-measure 1, this path exists (since $\mu(\mc{A}_{\smallNW})=1$), and has asymptotic slope $\sigma_{Nw}=-\frac{\mu(\mc{A}_{\smallN})}{1-\mu(\mc{A}_{\smallN})}$.  Similarly the $nW$ path from $x$, defined to be that path starting from $x$ obtained by following $\leftarrow$ whenever possible and otherwise following $\uparrow$, exists and has asymptotic slope $\sigma_{nW}=-\frac{1-\mu(\mc{A}_{\smallW})}{\mu(\mc{A}_{\smallW})}$. We can also define the $Se$ and $sE$ paths from $x$ and find their asymptotic slopes.  Except in the case that $\mu(\mc{A}_{\smallN})=1$ or $\mu(\mc{A}_{\smallS})=1$ (whence trivially $U$ or $L$ is infinite) we have that $-\infty<\sigma_{Nw}\le  \sigma_{nW}\le 0$ and $-\infty<\sigma_{Se}\le  \sigma_{sE}\le 0$.  Moreover, since every direction is possible we have that $\sigma_{Nw}<0$ and $\sigma_{Se}<0$.  

Let $k\ge 0$.  If $\sigma_{Nw}<\sigma_{sE}$, we may start from the origin and follow the $sE$ path until reaching a vertex $x_k$ on this path from which the $Nw$ path includes a vertex $(0,k+j)$ for some $j\ge 0$. This shows that $U_0$ is almost surely infinite.  If $\sigma_{nW}>\sigma_{Se}$ then follow the $Se$ path from the origin until reaching a vertex $y_k$ on this path from which the $nW$ path includes a vertex $(0,-k-j)$ for some $j\ge 0$. This shows that $L$ is almost surely infinite.  

It remains to verify the claim when $-\infty<\sigma_{Nw}=\sigma_{nW}=\sigma_{Se}=\sigma_{sE}<0$, which implies that $\mu(\mc{A}_{\smallN}\cap\mc{A}_{\smallW})=0$ and $\mu(\mc{A}_{\smallS}\cap\mc{A}_{\smallE})=0$.
In this case we define a new path, called the $nW_s$ path. From a vertex $x$, it follows whichever of $\leftarrow$ or $\uparrow$ is possible (now only one will be), except when at a vertex $y$ such that $\downarrow \in \mG_{y}$ and $\leftarrow \in \mc{G}_{y-e_2}$ (i.e.~such that $\mG_y\in \mc{A}_{\smallS}$ and $\mG_{y-e_2}\in \mc{A}_{\smallW}$). At such a vertex $y$, it follows the $\downarrow$ step, followed by the $\leftarrow$ step.  This path has a slope $\sigma_{nW_s}>\sigma_{Nw}=\sigma_{Se}$. We may now proceed as before by following the $Se$ path from $o$ until reaching a vertex $y_k$ from which the $nW_s$ path passes through $(0,-k-j)$ for some $j\ge 0$. This shows that $L$ is almost surely infinite, and completes the proof of part (a).

Part (b) now follows as in \cite[Proposition 3.8]{HS_DRE1}.
\qed\bigskip

The proof of Proposition \ref{prp:Ctrichotomy2} is a modification of that of \cite[Corollary 3.10]{HS_DRE1}, in exactly the same way as the proof of Proposition \ref{prp:Ctrichotomy} adapts that of \cite[Proposition 3.8]{HS_DRE1}.  We therefore omit it.

\section{Critical probabilities and coupling for Model \ref{exa:NE_SW}}
\label{sec:orthant_main}
For any site $x\in \Z^2$, let $\mathbf{C}_x$ denote the set of sites $y\in \Z^2$ for which there is some $N\ge 0$ and a sequence $\{x=y_0,y_1,y_2,\dots,y_N=y\}$ such that $\mc{G}_{y_i}=\NE$ for $0\le i\le N$ and $y_{i+1}-y_i\in \OTSP=\{-e_1,e_2,e_2-e_1\}$ for  $0\le i\le N-1$.  Similarly let $\mathbf{B}_x$ denote the set of sites $y\in \Z^2$ for which there is some $N\ge 0$ and a sequence $\{y=y_0,y_1,y_2,\dots,y_N=x\}$ such that $\mc{G}_{y_i}=\NE$ and $y_{i+1}-y_{i}\in \OTSP$ for each $i$.  Let 
\[\mathbf{A}=\{x\in \Z^2:|\mathbf{C}_x|=|\mathbf{B}_x|=\infty\},\]
i.e.~$\mathbf{A}$ is the set of sites that are in a bi-infinite cluster (in the sense of oriented site-percolation on the triangular lattice (OTSP)) of $\NE$ sites.  

For Model \ref{exa:NE_SW}, the conclusions of Proposition \ref{prp:Ctrichotomy} can be extended to the following.
\begin{PRP}
\label{prp:Cphasetransition}
For Model \ref{exa:NE_SW} $(\NE\SW)$:
\begin{enumerate}
\item $0\le p<1- p_c\quad \Rightarrow\quad \barc_o$ is a.s.~blocked above;
\item $1-p_c\le p\le p_c\quad\Rightarrow\quad \barc_o=\Z^2$ a.s.;
\item $p_c< p\le1\quad\Rightarrow\quad \barc_o$ is a.s.~blocked below.
\end{enumerate}
\end{PRP}
\proof We proceed as in the proof of \cite[Theorem 3.12]{HS_DRE1}, and will reiterate part of the latter in order to explain the role of the triangular lattice. Let $w(n)$ be decreasing. For $w$ to be a flbf for $\mc{G}$, the vertices in $w_\ge$ that have a (square-lattice) nearest-neighbour in $(w_{\ge})^c$ can be enumerated naturally as $\{y_t\}_{t \in \Z}$ to form a sequence of vertices moving upwards and to the left.  The possible transitions in this sequence of vertices are as follows. 
\begin{itemize}
\item Upwards, e.g.~from $(n,k)$ to $(n,k+1)$. This happens if $w(n)\le k<w(n-1)-1$. 
\item Leftwards, e.g.~from $(n,k)$ to $(n-1,k)$. This happens if $w(n)=w(n-1)=k$. 
\item Diagonally to the NW, e.g.~from $(n,k)$ to $(n-1,k+1)$. This happens if $w(n)\le k=w(n-1)-1$. 
\end{itemize}
These are three of the six possible transitions in a triangular lattice, whose families of lines are horizontal, vertical, and diagonal with slope $-1$ (the set of lattice points is still $\Z^2$). For $w(n)$ to be a flbf, it is necessary and sufficient that each vertex $y_t$ in this sequence has local environment $\NE$. Calling $\NE$ vertices ``open'' and $\SW$ vertices ``closed'', this sequence defines a bi-infinite oriented (triangular lattice) nearest-neighbour path such that $\mc{G}_{y_t}=\NE$.  It follows that each $y_t \in \mathbf{A}$ and in particular $(n,w(n))\in \mathbf{A}$.

Now continue exactly as in the proof of \cite[Theorem 3.12]{HS_DRE1} to obtain (c), substituting Proposition \ref{prp:Ctrichotomy} for \cite[Proposition 3.8]{HS_DRE1}.
A similar argument (or just symmetry) gives part (a), and part (b) follows immediately from the above argument and Proposition \ref{prp:Ctrichotomy}. 
\qed

\bigskip

For $p>p_c$ define 
$\del^*\mc{C}_o\subset \mc{C}_o$ by
\[\del^*\mc{C}_o =\{(n,m)\in \Z^2: m=W(n) \text{ or }W(n)\le m<W(n-1)\}.\]
The above sequence $y_t$ traces out $\del^*\mc{C}_o$ sequentially as a (triangular lattice) path. 
Clearly $o$ lies in or above the set $\del^*\mc{C}_o$, since $o \in \mc{C}_o$ implies that $W(0)\le 0$.  
We will need the following later. 
\begin{COR}
\label{cor:lbfcriterion}
For Model \ref{exa:NE_SW} $(\NE\SW)$, when $p>p_c$ the flbf $W(n)$ satisfies: 
\begin{align}
\label{Wsatisfies}
W(n)=\begin{cases}
K_0\equiv \sup\{k\le 0:(0,k)\in \mathbf{A}\} &, \text{ if }n=0,\\
\sup\{k\in \Z: (n,k)\in \mathbf{B}_{(0,K_0)}\cap \mathbf{A}\} &, \text{ if }n \in \N,\\
\sup\{k\in \Z: (n,k)\in \mathbf{C}_{(0,K_0)}\cap \mathbf{A}\} &, \text{ if }n \in -\N.
\end{cases}
\end{align}
\end{COR}
\proof Define $W^*(n)$ to be the right hand side of \eqref{Wsatisfies}.  Since $K_0$ is finite (as in the proof of Proposition \ref{prp:Cphasetransition}) and since if $x\in \mathbf{A}$ then $x+e,x-e'\in \mathbf{A}$ for some $e,e' \in \{e_1,-e_2,e_1-e_2\}$, it follows that $W^*:\Z\ra \Z$ is well defined. 
Our objective is to show that $W=W^*$. 

Note that $W^*$ is decreasing, since e.g.~if $n\ge 1$ and $W^*(n)=k$ then there exists $j\in [k,K_0]$ such that $(n-1,j)\in \mathbf{C}_{(n,k)}\cap \mathbf{B}_{(0,K_0)}$, and then $(n-1,j)\in \mathbf{B}_{(0,K_0)}\cap \mathbf{A}$.  Clearly then the origin cannot connect to anything below $W^*$, so $W^*\le W$.  We have seen in the proof of Proposition \ref{prp:Cphasetransition} that $(n,W(n))\in\mathbf{A}$ for every $n\in\Z$.  Since $W(0)\le 0$ and $(0,W(0))\in \mathbf{A}$, we get $W(0)\le W^*(0)$.  Therefore $W(0)=W^*(0)$. As in the proof of Proposition \ref{prp:Cphasetransition}, $(n,W(n))\in \mathbf{B}_{(0,W(0))}$ or $(n,W(n))\in \mathbf{C}_{(0,W(0))}$ (depending on the sign of $n\in\Z$) via the path $y_t$. Since $(n,W(n))\in\mathbf{A}$, we get that $W(n)\le W^*(n)$, and the result follows.
\qed

\bigskip 

We will need a corresponding result for an infinite $\mc{B}_o$ cluster.  Indeed, by \cite[Theorem 3.12]{HS_DRE1} there exists a  decreasing function $V:\Z\ra \Z$ with $\mc{B}_o=V_<$ (in \cite{HS_DRE1}, $V-1$ is a bubf, so $V$ is a flbf).  
The proof of the following proceeds just as in that of Corollary \ref{cor:lbfcriterion}, so will be omitted. 
\begin{COR}
\label{cor:lbfcriterionforB}
Consider Model \ref{exa:NE_SW} $(\NE\SW)$ with $p>p_c$. If $\mc{B}_o$ is infinite, then $\mc{B}_o=V_<$ for a decreasing flbf $V:\Z\to\Z$ satisfying:
\begin{align*}
V(n)=&\begin{cases}
K'_0\equiv \inf\{k\ge 0:(0,k)\in \mathbf{A}\}, & \text{ if }n=0,\\
\inf\{k\in \Z: (n,k)\in \mathbf{B}_{(0,K'_0)}\cap \mathbf{A}\}, & \text{ if }n \in \N,\\
\inf\{k\in \Z: (n,k)\in \mathbf{C}_{(0,K'_0)}\cap \mathbf{A}\}, & \text{ if }n \in -\N.
\end{cases}
\end{align*}
\end{COR}

As $p$ increases from 0 to 1, simulation suggests that the cluster $\mc{B}_o$ can change from infinite to finite and back again arbitrarily many times, while ``holes" in $\mc{C}_o$ can expand and contract. Nevertheless, both clusters have a kind of monotonicity, as in Theorem \ref{thm:orthant_monotone}. \bigskip
 
\noindent {\em Proof of Theorem \ref{thm:orthant_monotone}.} Couple the environments for all $p$ as in (\ref{eq:coupling}) so that as $p$ decreases we switch $\NE$ to $\SW$.  Due to Proposition \ref{prp:Cphasetransition}, the conclusion of the Theorem is trivial if $p_1=p_c$.  So assume $p_c<p_1<p_2\le 1$, and let $W$ be the flbf with $\mc{C}_o(p_1)=W_\ge$. This implies that every site  $x\in \del^*\mc{C}_o(p_1)$ satisfies $x\in \mc{C}_o(p_1)$ and $\mc{G}_x(p_1)=\NE$. Therefore $\mc{G}_x(p_2)=\NE$ for all such $x$, which implies that regardless of $\{\mc{G}_y(p_2): y \notin \del^*\mc{C}_o(p_1)\}$, there can be no open path in $\mc{G}(p_2)$ from $o$ to any site below $\del^*\mc{C}_o(p_1)$.  In other words, $W_<\subset (\mc{C}_o(p_2))^c$, and by definition of $\barc_o$ also $W_<\subset (\barc_o(p_2))^c$.  This implies that $(\barc_o(p_1))^c\subseteq (\barc_o(p_2))^c$, which establishes the desired conclusion except for showing that the inequality is strict.  

It remains to prove strictness, ie. $\barc_o(p_2)\neq \barc_o(p_1)$.  Let $G_x=\{U_x\le p_2\}=\{\mc{G}_x(p_2)=\NE\}$.  Let $\mc{H}_2=\sigma(\{G_x,x \in \Z^2\})$.  Then $\mc{C}_o(p_2)$ and $\del^*\mc{C}_o(p_2)$ are $\mc{H}_2$-measureable random sets, and $\del^*\mc{C}_o(p_2)\subset H_2=\{x\in \Z^2: \mc{G}_x(p_2)=\NE\}$.  Conditional on $\mc{H}_2$, for any $\mc{H}_2$-measureable subset $H'$ of $H_2$ we have that $\{U_x:x\in H'\}$ are i.i.d.~$U[0,p_2]$ random variables. In particular $\{U_x:x\in \del^*\mc{C}_o(p_2)\}$ are i.i.d.~$U[0,p_2]$ random variables under this conditioning.  Thus, $\mP\big(\{x\in \del^*\mc{C}_o(p_2):U_x>p_1\}=\varnothing\big|\mc{H}_2\big)=0$ almost surely, so $\mP\big(\{x\in \del^*\mc{C}_o(p_2),U_x>p_1\}=\varnothing\big)=0$.  This says that (almost surely) there exists $u\in \del^*\mc{C}_o(p_2)$ with $U_u\in (p_1,p_2]$.  Therefore $\mc{G}_u(p_1)=\SW$, so $u\notin \del^*\mc{C}_o(p_1)$, hence $\del^*\mc{C}_o(p_1)\ne \del^*\mc{C}_o(p_2)$.  
\Qed

The following is a version of Theorem \ref{thm:orthant_monotone} for the clusters $\mc{B}_o(p)$. The proof is similar, and will be omitted.
\begin{COR}
\label{cor:orthantB_monotone}
Consider Model \ref{exa:NE_SW} $(\NE\SW)$, with the coupling (\ref{eq:coupling}). Let $p_c\le p_1<p_2\le 1$. Then $\mc{B}_o(p_1)\supsetneq \mc{B}_o(p_2)$ $\mP$-a.s. on $\{|\mc{B}_o(p_1)|=\infty\}$.
\end{COR}

\section{Oriented site percolation on the triangular lattice}
\label{sec:OP}
In this section we state without proof a number of results about the OTSP model $(\OTSP,\cdot)$ that follow using the methods of  \cite{Dur84} for two dimensional oriented percolation models. An expanded version of this paper, with the omitted proofs included, is available from the authors upon request.  In this model we have local environment $\mathbf{G}_x=\OTSP$ with probability $p$, and $\mathbf{G}_x=\emptyset$ with probability $1-p$, both on the triangular lattice described in Section \ref{sec:orthant_main}. Recall that forward clusters in this model are denoted $\mathbf{C}_x$, and backward clusters $\mathbf{B}_x$. The natural coupling (\ref{eq:coupling}) gives a probability space on which the sets $\mathbf{C}_o(p)$ are increasing in $p$ almost surely, so
\[\Theta_+(p)=\mP(|\mathbf{C}_o(p)|=\infty) \quad \text{ is increasing in }p,\]
giving the critical value $p_c^{\smallOTSP}=\inf\{p:\Theta_+(p)>0\}\in (0,1)$.  

The estimated value is $p_c^{\smallOTSP}\approx 0.5956$ (see De Bell and Essam \cite{DeBE} or Jensen and Guttmann \cite{JG}). 
The best rigorous bounds that we have found in the literature are $0.5466\le p_c^{\smallOTSP}\le 0.7491$, where the latter comes from the fact that $p_c^{\smallOTSP}\le p_c^{\smallNE}\le 0.7491$ (the latter, referring to oriented percolation on the square lattice, is in Balister et al \cite{BBS}). 
Similarly, if $p_c^{\text{TSP}}$ denotes the critical threshold for (un-oriented) triangular site percolation, then 
$p_c^{\smallOTSP}\ge p_c^{\text{TSP}}=1/2$ (see Hughes \cite{Hughes}). The lower bound of $0.5466$ comes from \cite{HS_DRE1} based on estimates of the critical value for the model $(\NE\,\SW)$.  
In Section \ref{sec:pcbounds} we improve this lower bound to $.5731\le p_c^{\smallOTSP}$, by finding arbitrarily large connected circuits in the dual model $(\NE\,\SW)$ when $0.4269\le p\le 0.5731$

In order to describe the shape of an infinite $\mathbf{C}_x$ cluster, define $w_n=\sup\{x:(-n,x)\in \mathbf{C}_o\}$ and $v_n=\inf\{x:(-n,x)\in \mathbf{C}_o\}$.  The following Proposition is proved using subadditivity of quantities related to $w_n$. Minor modifications arise from the proofs in \cite{Dur84}, because the latter treats oriented bond percolation on the square lattice, while we need oriented site percolation on the triangular lattice.
\begin{PRP}
\label{prp:OPslope}
For the percolation model $(\OTSP,\cdot)$ with $1>p>p_c^{\smallOTSP}$, there exists $\rho=\rho_p<-1$ such that almost surely on the event $\{|\mathbf{C}_o|=\infty\}$, the upper and lower boundaries of $\mathbf{C}_o$ have asymptotic slopes $\rho$ and $1/\rho$ respectively. In other words, $\frac{w_n}{-n}\ra \rho$ and $\frac{v_n}{-n}\ra 1/\rho$ almost surely as $n\ra \infty$.  
\end{PRP}

Since $v_n$ is bounded below by a sum of independent Geometric$(1-p)$ random variables, we get the inequality $-\frac{p}{1-p}\le \rho_p$.  The following two additional Lemmas can be proved as in \cite{Dur84}.
\begin{LEM}
\label{lem:OPslope}
$\rho_p$ is continuous and strictly decreasing in $p> p_c^{\smallOTSP}$, with $\rho_p\uparrow -1$ as $p\downarrow p_c^{\smallOTSP}$.
\end{LEM}
Let $\tau=\sup\{y-x: (x,y)\in \mathbf{C}_o\}$, which measures the furthest diagonal line reached by the forward cluster of the origin.  More generally, if $z=(x_0,y_0)$, let 
$\tau_z=\sup\{(y-y_0)-(x-x_0): (x,y)\in \mathbf{C}_z\}$. Note that $|\mathbf{C}_o|=\infty\Leftrightarrow \tau=\infty$.
\begin{LEM}
\label{lem:dual_tail}
If $p>p_c^{\smallOTSP}$, then there exist constants $C$, $\gamma>0$ such that $\mP(n\le \tau<\infty)\le Ce^{-\gamma n}$.
\end{LEM}

On the event $\{|\mathbf{C}_o|=\infty\}$, let $a_n=\sup\{x_n^{[2]}:x_n\in \mathbf{C}_o, |\mathbf{C}_{x_n}|=\infty, x_n^{[1]}=-n\}$ for $n \in \Z_+$.  The following result says that $a_n$ has the same asymptotic slope as the upper boundary of $\mathbf{C}_o$.
\begin{COR}
\label{cor:a_op_slope}
For the percolation model $(\OTSP,\cdot)$ with $p>p_c^{\smallOTSP}$, $\lim \frac{a_n}{-n}=\rho_p$ almost surely on the event $\{|\mathbf{C}_o|=\infty\}$.
\end{COR}
\proof 
Let $p>p_c^{\smallOTSP}$.
Since $a_n\le w_n$ for every $n$, and $\frac{w_n}{-n}\to\rho_p>-\infty$, 
it suffices to prove that for each $\epsilon>0$, 
\[\mP(w_n-a_n>\epsilon n \, \text{ infinitely often},    |\mathbf{C}_o|=\infty)=0.\]
Let $R>-\rho_p$. Then $w_n\le Rn$ for all sufficiently large $n$, so we may find $N$ such that $w_n\le Rn$ for every $n\ge N$. It will therefore suffice to show that 
\begin{equation}
\label{eqn:iocondition}
\mP(\text{$|\mathbf{C}_o|=\infty$, and $w_n-a_n>\epsilon n$ for infinitely many $n\ge N(1+R)$} )=0.
\end{equation}
Call $y-x$ the {\em generation} of a point $(x,y)\in\mc{L}$. Along any open path in this model, the generation increases at each step, by 1 or by 2. Fix $\epsilon>0$, and $n>4/\epsilon$. Suppose that $|\mathbf{C}_o|=\infty$, and $N(1+R)\le n$, and $w_n-a_n> \epsilon n$. Let $T$ be the triangle with vertices $a=(-n,0)$, $b=(-n,Rn)$, and $c=(-N,RN)$, and let $T_0\supset T$ be the triangle with vertices $a$, $b$, and $o$. Any  open path $\Gamma$ from $o$ to $(-n,w_n)$ enters $T$ along the side $ac$. Since $n\ge N(1+R)$, $ac$ has slope $\le 1$, so entry to $T$ occurs no later than generation $n$. Therefore from generation $n$ through $n+w_n$ the path $\Gamma$ lies entirely within $T$. Consider the lattice point $z=(x,y)$ on this path whose generation first exceeds $n+w_n-\epsilon n$. Then $w_n-\epsilon n>a_n\ge 0$, so $n+w_n-\epsilon n>n$, and hence $z\in T$. The generation of $(-n,w_n)$ is at least $\epsilon n -2>\epsilon n/2$ larger than that of $z$, so $\tau_z\ge\epsilon n/2$. On the other hand, $x\ge -n$, so $y=(y-x)+x\ge n+ w_n-\epsilon n+x\ge w_n-\epsilon n>a_n$. Therefore we cannot have $|\mathbf{C}_z|=\infty$. In other words, $\tau_z<\infty$. We have shown that for $n>4/\epsilon$,
$$
\{\text{$|\mathbf{C}_o|=\infty$, $N(1+R)\le n$, and $w_n-a_n> \epsilon n$}\}
\subset\bigcup_{z\in T_0}\Big\{\frac{\epsilon n}{2}\le\tau_z<\infty\Big\}.
$$
There are at most $Rn^2$ lattice points in $T_0$, so by Lemma \ref{lem:dual_tail}, the probability of this event is at most $CRn^2e^{-\gamma n\epsilon/2}$, which sums. Therefore \eqref{eqn:iocondition} follows by Borel-Cantelli.\qed

\section{Asymptotic slopes for Model \ref{exa:NE_SW}}
\label{sec:slopes}
To complete this circle of results, it simply remains to show that the remaining parts of Theorem \ref{thm:orthant_main} follows from the results of the previous section.

\bigskip

{\noindent \em Proof of Theorem \ref{thm:orthant_slope}.}
By Proposition \ref{prp:Cphasetransition}, when $p>p_c$, $\barc_o$ is bounded below by a flbf $W(n)$. As in Corollary \ref{cor:lbfcriterion}, we construct $y\equiv (0,K_0)\in {\bf A}$.  So $\{|{\bf C}_{y}|=\infty=|{\bf B}_{y}|\}$ and for $n \in \N$, we may define $a_n(y)=\sup\{z_n^{[2]}:z_n\in {\bf C}_{y},|{\bf C}_{z_n}|=\infty, z_n^{[1]}=-n\}$.   Corollary \ref{cor:lbfcriterion}  now implies that $W(-n)=a_n(y)$, so by Corollary \ref{cor:a_op_slope} and translation invariance, we get that $\frac{W(-n)}{-n}\ra \rho_p$, where $\rho_p$ is as in Proposition \ref{prp:OPslope}.  By that result and Lemma \ref{lem:OPslope}, $\rho_p\uparrow -1$ as $p \downarrow p_c$. This establishes the desired statements for the northwest boundary. The results for the southeast boundary follow by symmetry. \qed

\bigskip

{\noindent \em Proof of Theorem \ref{thm:Mphasetransition}.}  By \cite[Theorem 4.9(b)]{HS_DRE1}, it simply remains to show that $\mc{M}_o$ is a.s.~finite when $p>p_c$.  Let $p>p_c$.  Then as above, $\barc_o$ is bounded below by $W$ which satisfies $\frac{W(-n)}{-n} \ra \rho<-1$ and by symmetry in the model $(\NE,\SW)$ also $\frac{W(n)}{n} \ra 1/\rho>-1$ for $n \in \N$.  Similarly by Corollary \ref{cor:lbfcriterionforB}, $\mc{B}_o$ is bounded above by $V$, and $\frac{V(-n)}{-n} \ra 1/\rho$ by symmetry in the dual oriented site percolation model. Then also $V(n)/n\ra \rho$ by symmetry in the model $(\NE,\SW)$.  Since $\rho<-1$ it follows that $V(n)< W(n)$ for all but finitely many $n$, so $\mc{M}_o$ is finite. See Figure \ref{fig:orthant_C_o_B_o}.\qed

\section{Lower bounds on $p_c$}
\label{sec:pcbounds}

\bigskip

We define a nearest-neighbour (self-avoiding) path $\vec{x}=x_0,x_1,x_2,\dots$ to be a {\em $\C_x$-path} if $x_0=x$ and $x_{i+1}-x_i\in \mc{G}_{x_i}$ for each $i\in \Z_+$.  For $\sigma\in \re$ and $x\in \Z^2$ we write $\mc{S}_{x}(\sigma)$ for the set of $\C_x$-paths $\vec{x}$ such that $x_n^{[2]}\ra \infty$ and $x_n^{[2]}/x_n^{[1]}\ra \sigma$ as $n\ra \infty$.  Let $\mc{S}(\sigma)=\mc{S}_{o}(\sigma)$.  Similarly define  $\mc{S}^-_{x}(\sigma)$ [resp. $\mc{S}^+_{x}(\sigma)$] as the set of $\C_x$-paths  such that $x_n^{[2]}/x_n^{[1]}\ra \sigma$ and $x_n^{[1]}\ra -\infty$ [resp. $+\infty$]. Write $\mc{S}^-(\sigma)=\mc{S}^-_{o}(\sigma)$ and $\mc{S}^+(\sigma)=\mc{S}^+_{o}(\sigma)$.

The following Lemma gives a strategy for generating better one-sided bounds for the critical point for the orthant model, and hence also for oriented site percolation on the triangular lattice. 
\begin{LEM}
\label{lem:lookforslopes}
Consider the model $(\NE \SW)$, with $p\ge  .5$.  Suppose that for some $\sigma>-1$, there exists a path $\vec{x}\in \mc{S}^-(\sigma)$ almost surely.  Then $p\le p_c$.
\end{LEM}
\proof Suppose that $p>p_c$.  Then $\barc_o$ is bounded below by $W$, where $\frac{W(-n)}{-n}\to\rho_p<-1$ as in Proposition \ref{prp:OPslope}.   It follows that for any $\sigma>\rho_p$, the set $\{x \in \mc{C}_o: x_i^{[1]}< 0,x_i^{[2]}/x_i^{[1]}>\sigma\}$ is finite.  Since $\rho_p<-1$, this implies that there can be no infinite path $\vec{x}\in \mc{S}^-(\sigma)$ for $\sigma>-1$.\qed

Let
\[g_1(p)=\frac{p^2}{1-p}-\frac{1-p}{p}=\frac{p^3-(1-p)^2}{p(1-p)}=\frac{p^3-p^2+2p-1}{p(1-p)}.\]

\begin{LEM}
\label{lem:Cslope}
For the model $(\NE,\leftarrow)$, $\mc{C}_o$ contains (self-avoiding) $\C_o$-paths $L_o\in \mc{S}(\frac{-p}{1-p})$ and $R_o\in \mc{S}(g_1(p)^{-1})$.  Moreover, environment-connected components of  $\mc{C}_o^c$ between these two paths are finite.
\end{LEM}
\proof  Define the {\em NW path} from $x$ to be the path obtained by always choosing to follow $\leftarrow$ when possible and otherwise following $\uparrow$.  
The asymptotic slope of the NW path (call this path $L_o$) from the origin is $-\frac{p}{1-p}$, which establishes the first claim. 
For the second claim, consider the path $R_o$ from the origin that evolves as follows.  Whenever the environment at the current location is $\leftarrow$, the path follows this west arrow.  Whenever the environment at the current location $x$ is $\NE$, the path follows the east arrow to $x+(1,0)$ if the environment at $x+(1,0)$ is also $\NE$, otherwise the path follows the north arrow.  By definition this path never backtracks, so it is self-avoiding.  After each northern step taken by the path, the environment thereafter encountered has never been viewed before, hence each northern step constitutes a renewal.  Thus the path moves upwards through the set of horizontal bands $\Z\times\{k\}$, $k=0,1,2,\dots$.  Let $(X_k,k)$ be the point where our path first enters the $k$th band. Then we can represent $X_{k+1}$ as follows:
\begin{multline*}
X_{k+1}=\sup\{j:\text{$\exists$ path from $(X_k,k)$ to $(j,k+1)$ consistent with the environment} \\
\text{and lying within the $k$th band, except for the final step}\}.
\end{multline*}
Then the $\Delta_k=X_{k+1}-X_k$ are i.i.d.~with
\[\mP(\Delta_0=k)=\begin{cases}
p^{k+1}(1-p), & \text{ if }k\ge 0\\
(1-p)^kp, & \text{ if }k<0.
\end{cases}\] 
It follows that 
\[\mE[\Delta_0]=\sum_{k=0}^\infty kp^{k+1}(1-p)-\sum_{k=1}^{\infty}k(1-p)^kp=\frac{p^2}{1-p}-\frac{1-p}{p}=g_1(p),\]
and that the asymptotic slope of the path $R_o$ is $1/g_1(p)$ as claimed. 
 
By construction, a vertex $(x,n)\in R_o$ never lies strictly to the left of any vertex $(j,n)\in L_o$ (the paths $R_o$ and $L_o$ may meet, but not cross).  An elementary comparison of the asymptotic slopes shows that $R_o=\{r_0,r_1,\dots\}$ eventually lies strictly to the right of $L_o=\{l_0,l_1,\dots\}$, in the sense that there exists some $m$ such that for all $n\ge m$, $\inf\{j:(j,n)\in R_o\}>\sup\{j:(j,n)\in L_o\}$.  Trivially also by construction $r_n^{[2]}\ra \infty$ and similarly for $l_n^{[2]}$.  To prove the last claim observe that the NW path from each $r_n$ eventually hits the path $L_o$ (since NW paths from any two vertices eventually meet), and that these NW paths are all in $\mc{C}_o$.   \qed
 
Note that it follows immediately from this result that we can find similar paths (up to symmetry) in $\mc{C}_o$ for the models containing this as a submodel.  This implies results for other models, such as the following improvement on the lower bound on $p_c$ in \cite[Theorem 4.12]{HS_DRE1}.
\begin{COR}
\label{cor:Corthant}
In Model \ref{exa:NE_SW} $(\NE \SW)$, $\barc_o=\Z^2$ for $p \in [.5,.5699)$. Therefore $p_c^{\smallOTSP}\ge .5699$.
\end{COR}
\proof Since Model \ref{exa:NE_SW} contains $(\uparrow \SW )$, by Lemma \ref{lem:Cslope} and symmetry, $\mc{C}_o$ contains a self-avoiding path $P_1\in \mc{S}^-(g_1(1-p))$.  So $\barc_o=\Z^2$ for any $p\ge 0.5$ such that $g_1(1-p)>-1$, by Proposition \ref{prp:Cphasetransition} and Lemma \ref{lem:lookforslopes}. 
The condition $g_1(1-p)>-1$ holds for $p \in [\hlf,p_1)$, where $p_1^3-p_1^2+2p_1-1=0$, i.e.~for $p \in [\hlf, .5699)$ 
(the condition is actually equivalent to $g_1(p)<0$ since $g_1(p)+g_1(1-p)=-1$).  
\qed

\medskip

The true value we are aiming for is $p_c^{\smallOTSP}$, which is estimated as $0.5956$. The path $R_o$ above is defined (see Lemma \ref{lem:Cslope}) in terms of expected horizontal displacements that can be achieved for paths consistent with the environment, staying within a horizontal band of width 1. 
In fact, there is a sequence of estimates that in principle should converge to the true value. Repeat the above argument, but using bands of width $K$. The resulting bounds $p_K$ should converge to $p_c^{\smallOTSP}$.  We'll content ourselves with computing $p_2$. 

\bigskip

{\noindent \em Proof of Theorem \ref{thm:OTSP_pc}.}  
In analogy with the proof of Lemma \ref{lem:Cslope}, we let $X_0=0$ and $\Delta_k=X_{k+1}-X_k$, where
\begin{multline*}
X_{k+1}=\sup\{j:\text{$\exists$ path from $(X_k,2k)$ to $(j,2k+2)$ consistent with the environment} \\
\text{and lying within the band $\Z\times\{2k,2k+1\}$, except for the final step}\}.
\end{multline*}
We set $g_2(p)=\mE[\Delta_0]/2$ so that the path using horizontal bands has slope $1/g_2(p)$. The strategy of Corollary \ref{cor:Corthant} works with $g_2$ in place of $g_1$, provided $g_2(1-p)>-1$.  Therefore we are left to compute $g_2(p)$.  This is done below, ultimately leading to $p_2=1-q$ where
\begin{align}
{q}^{11}-6\,{q}^{10}+18\,{q}^{9}-38\,{q}^{8}+64\,{q}^{7}-90\,{q}^{6}+104\,{q}^{5}-94\,{q}^{4}+66\,{q}^{3}-34\,{q}^{2}+12\,q-2=0.\label{polynom}
\end{align}
Using Newton's method (with a starting point of $q=.4$, in the computer algebra system {\em Maxima}) gives $p_2\approx .5730$, which completes the proof.  In the remainder of this section, we show how to obtain \eqref{polynom}.

We write $\mc{G}_0=\binom{\mc{G}_{(0,1)}}{\mc{G}_{(0,0)}}$ for the pair of local environments at $(0,0)$ and $(0,1)$. First consider the case $\mc{G}_0=\binom{\smallNE}{\smallNE}$. Our first goal will be to compute
$$
\mu_{\smallNE}^{\smallNE}=\mE\Big[\Delta_0\mid \mc{G}_0=\binom{\smallNE}{\smallNE}\Big].
$$
Paths are not actually unique, but we take the convention that in this situation the path moves $\uparrow$ from $(0,0)$ to $(0,1)$, and then considers its next move. If $\mc{G}_1=\binom{\smallNE}{\cdot}$  (the $\cdot$ simply means an arbitrary environment) the path moves $\rightarrow$ to $(1,1)$. If $\mc{G}_1=\binom{\smallSW}{\smallSW}$ then we have hit an impassable obstacle, and the path has no choice but to exit using $\uparrow$, in which case $X_1=0$. The remaining possibility is $\mc{G}_1=\binom{\hspace{-.2cm}\smallSW}{\hspace{.2cm}\smallNE}$. There could in fact be a sequence of such pairs, followed either by a $\binom{\cdot}{\smallSW}$ or by $\binom{\smallNE}{\smallNE}$. For example, 
$$
\mc{G}_0\mc{G}_1\dots \mc{G}_5=
\begin{matrix}
\smallNE & \ssmallSW & \ssmallSW & \ssmallSW & \ssmallSW & \cdot \\
\smallNE & \smallNE & \smallNE & \smallNE & \smallNE & \hspace{-.2cm}\smallSW
\end{matrix}
$$
or 
$$
\mc{G}_0\mc{G}_1\dots \mc{G}_5=
\begin{matrix}
\smallNE & \ssmallSW & \ssmallSW & \ssmallSW & \ssmallSW & \smallNE \\
\smallNE & \smallNE & \smallNE & \smallNE & \smallNE & \smallNE
\end{matrix}
$$
The first case also represents an impassable obstacle, so our path chooses to exit using $\uparrow$, making $X_1=0$, In the second case, we move $\rightarrow$, $\downarrow$, take a sequence of $\rightarrow$'s, and then move $\uparrow$ to reach the top row of the band again. 

This description makes it clear that there is a renewal structure here, with the construction starting afresh every time we reach a new $\NE$ environment on the top row of the band. To formalize this, let $(Z_j,1)$ be the site of the $j$th such new $\NE$ environment reached by our path (with $Z_0=0$ corresponding to the initial environment). If $J$ denotes the total number of such environments reached, then we take $Z_j=Z_J$, $\forall j\ge J$. The dynamics are that 
$$
Z_{j+1}=
\begin{cases}
Z_{j}+1,&\text{if $\mc{G}_{Z_j+1}=\binom{\smallNE}{\cdot}$}\\
Z_{j}+i,&\text{if } \mc{G}_{Z_j+1}\dots \mc{G}_{Z_j+i}=
\begin{matrix}
\ssmallSW & \dots & \ssmallSW & \smallNE\\
\smallNE & \dots & \smallNE & \smallNE
\end{matrix}
\\
Z_{j}, &\text{if }\mc{G}_{Z_j+1}=\binom{\smallSW}{\smallSW} \text{ or } \mc{G}_{Z_j+1}\ldots = 
\begin{matrix}
\ssmallSW & \dots & \ssmallSW & \cdot\\
\smallNE & \dots & \smallNE & \hspace{-.2cm}\smallSW
\end{matrix}
\end{cases}
$$
Let $\eta=p(1-p)$ and
$$
\alpha=\frac{p^2}{p^2+1-p}=\frac{p^2}{1-p(1-p)}=\frac{p^2}{1-\eta}
$$
denote the probability of encountering a $\binom{\smallNE}{\smallNE}$ before a $\binom{\cdot}{\smallSW}$. Let $\theta=p+p(1-p)\alpha$ denote the probability that $Z_1>Z_0$. Then 
$$
\mE[Z_1]=p+p(1-p)\sum_{i=1}^\infty (i+1)[p(1-p)]^{i-1}p^2=\theta+\frac{\alpha p(1-p)}{1-p(1-p)}
$$
and the renewal structure implies that
$$
\mu_{\smallNE}^{\smallNE}
=\mE\Big[\sum_{j=0}^\infty [Z_{j+1}-Z_{j}]1_{j\le J}\Big]=\sum_{j=0}^\infty \mE[Z_1]\theta^j=\frac{\mE[Z_1]}{1-\theta}.
$$
The second case we consider is that $\mc{G}_0=\binom{\hspace{-.2cm}\smallSW}{\hspace{.2cm}\smallNE}$, so our next goal is to compute
$$
\mu_{\smallNE}^{\smallSW}=\mE\Big[\Delta_0\mid \mc{G}_0=\binom{\hspace{-.2cm}\smallSW}{\hspace{.2cm}\smallNE}\Big].
$$
 If 
 $$
 \mc{G}_{0}\ldots \mc{G}_j =
\begin{matrix}
\ssmallSW & \dots & \ssmallSW & \smallNE\\
\smallNE & \dots & \smallNE & \smallNE
\end{matrix}
$$
then we travel $j$ steps $\rightarrow$ and then $\uparrow$ and find ourselves back in the situation just considered. While if 
$$
\mc{G}_0\mc{G}_1\ldots = 
\begin{matrix}
\ssmallSW & \dots & \ssmallSW & \cdot\\
\smallNE & \dots & \smallNE & \hspace{-.2cm}\smallSW
\end{matrix}
$$
then we are blocked to the right, and instead travel $\uparrow$ and then follow $\leftarrow$'s till reaching a $\NE$ on the top row of the band, at which point we exit from the band via that $\uparrow$. This gives us the expression
$$
\mu_{\smallNE}^{\smallSW}=\sum_{j=1}^\infty j[p(1-p)]^{j-1}p^2 +\alpha \mu_{\smallNE}^{\smallNE}-(1-\alpha)\sum_{i=1}^\infty i(1-p)^{i-1}p
=\alpha \mu_{\smallNE}^{\smallNE}+\frac{\alpha}{1-p(1-p)}-\frac{1-\alpha}{p}.
$$

The third case of interest is 
$$
\mu_{\smallSW}^{\smallSW}=\mE\Big[\Delta_0\mid \mc{G}_0=\binom{\smallSW}{\smallSW}\Big].
$$
We start out by following $\leftarrow$ along the bottom row, till reaching a $\NE$ site, at which point we go $\uparrow$ to the top row. If the site so reached is a $\SW$ then we again proceed $\leftarrow$ till reaching a $\NE$ on the top row, at which point we exit from the band via the $\uparrow$. On the other hand, if the first site reached in the top row is a $\NE$ then we have the opportunity to regain some lost ground. We step $\rightarrow$ along the top row as long as possible, and only go $\uparrow$ just before reaching a $\SW$. 
For example, 
$$
\mc{G}_{-5}\mc{G}_{-4}\mc{G}_{-3}\mc{G}_{-2} \mc{G}_{-1}\mc{G}_0=
\begin{matrix}
\smallNE & \smallNE & \smallNE & \smallNE & \hspace{-.2cm} \smallSW & \smallSW \\
\smallNE & \ssmallSW & \ssmallSW & \ssmallSW &\hspace{-.2cm} \smallSW & \smallSW
\end{matrix}
$$
has $X_1=-2$. 
This leads to the general expression
\begin{align*}
\mu_{\smallSW}^{\smallSW}&=
\sum_{k=1}^\infty (1-p)^{k-1}p\Big[ -k-\sum_{j=1}^\infty j(1-p)^jp+\sum_{j=0}^{k-2}jp^{j+1}(1-p)+(k-1)p^k\Big]\\
&=-\frac{1}{p}-\frac{1-p}{p}+\sum_{j=0}^\infty jp^{j+2}\sum_{k=j+2}^\infty (1-p)^k +\sum_{k=1}^\infty (k-1)p^2[p(1-p)]^{k-1}\\
&=-\frac{2-p}{p}+\frac{\alpha(1-p)}{1-p(1-p)}[p+(1-p)^2].
\end{align*}

The fourth and final case of interest is 
$$
\mu_{\smallSW}^{\smallNE}=\mE\Big[\Delta_0\mid \mc{G}_0=\binom{\,\,\smallNE}{\,\,\ssmallSW}\Big].
$$
Once again, we must go $\leftarrow$ till reaching the first $\NE$ on the bottom row of the band, and then go $\uparrow$. If this leads to a $\SW$ vertex, then all we can do is head $\leftarrow$ on the top row, till reaching a $\NE$ vertex, at which point we may exit via $\uparrow$. However, if we reach the top row at a $\NE$ vertex then we can regain lost ground by heading $\rightarrow$. Either this ends as in the previous case, before making it all the way back to $(0,1)$. Or we do make it back to $(0,1)$ this way, in which case we find ourselves back in the first case examined above. For example, 
$$
\mc{G}_{-5}\mc{G}_{-4}\mc{G}_{-3}\mc{G}_{-2} \mc{G}_{-1}\mc{G}_0=
\begin{matrix}
\smallNE & \smallNE & \smallNE & \smallNE & \smallNE & \smallNE \\
\smallNE & \ssmallSW & \ssmallSW & \ssmallSW & \ssmallSW & \ssmallSW
\end{matrix}
$$
has $X_1\ge 0$. This leads to the general expression
\begin{align*}
\mu_{\smallSW}^{\smallNE}&=
\sum_{k=1}^\infty (1-p)^{k-1}p\Big[ -k-\sum_{j=1}^\infty j(1-p)^jp+\sum_{j=0}^{k-2}jp^{j+1}(1-p)+\big(k+\mu_{\smallNE}^{\smallNE}\big)p^k\Big]\\
&=-\frac{2-p}{p}+\frac{\alpha}{1-p(1-p)}[1+(1-p)^3]+\alpha \mu_{\smallNE}^{\smallNE}.
\end{align*}
Thus, $2g_2(p)$ is equal to 
\begin{align}
\mE[\Delta_0]=&p^2\mu_{\smallNE}^{\smallNE} +p(1-p)\mu_{\smallNE}^{\smallSW} +p(1-p)\mu_{\smallSW}^{\smallNE} +(1-p)^2\mu_{\smallSW}^{\smallSW}.
\end{align}
This gives us that 
\[g_2(p)=\frac{{p}^{11}-6\,{p}^{10}+16\,{p}^{9}-30\,{p}^{8}+46\,{p}^{7}-62\,{p}^{6}+72\,{p}^{5}-66\,{p}^{4}+48\,{p}^{3}-26\,{p}^{2}+10\,p-2}{2\,{p}^{9}-8\,{p}^{8}+18\,{p}^{7}-28\,{p}^{6}+32\,{p}^{5}-28\,{p}^{4}+18\,{p}^{3}-8\,{p}^{2}+2\,p}.\]
Here the denominator is equal to 
\[2\,{\left( p-1\right) }^{2}\,p\,\left( {p}^{2}+1\right) \,{\left( {p}^{2}-p+1\right) }^{2}.\]
Solving for $g_2(1-p)=-1$ is equivalent to finding $p=1-q$, where 
\[{q}^{11}-6\,{q}^{10}+18\,{q}^{9}-38\,{q}^{8}+64\,{q}^{7}-90\,{q}^{6}+104\,{q}^{5}-94\,{q}^{4}+66\,{q}^{3}-34\,{q}^{2}+12\,q-2=0,\]
as claimed.
\qed

 \section*{Acknowledgements}
Holmes's research is supported in part by the Marsden fund, administered by RSNZ. Salisbury's research is supported in part by NSERC. Both authors acknowledge the hospitality of the Fields Institute, where part of this research was conducted.

\bibliographystyle{plain}

\end{document}